\documentclass{amsart}
\usepackage{ObermanLatexPackage}
\usepackage{cleveref}
\usepackage{mathtools}

\usepackage{siunitx} 
\newcommand{\Td}{{\mathbb{T}^d}}
\newcommand{\T}{{\mathbb{T}}}
\newcommand{\Sd}{{\mathcal{S}^d}}
\newcommand{\Ac}{{\mathcal{A}}}
\newcommand{\PrAB}{\Pr(\Td \times \Ac)}
\newcommand{\al}{\alpha}
\newcommand{\dd}{\,\mathrm{d}}
\newcommand{\HM}[0]{\operatorname{HM}}
\newcommand{\Hb}{\overline{H}(Q)}

\newcommand{\lp}{\lambda_\text{max}}
\newcommand{\lm}{\lambda_\text{min}}

\begin{document}
\title[Homogenization of Fully Nonlinear PDEs]{Approximate homogenization of convex nonlinear elliptic PDEs}
\author{Chris Finlay}
\author{Adam~M. Oberman}
\date{\today}
\thanks{The authors are thankful to Panagiotis Souganidis,  Guillaume Carlier, and Yifeng Yu for valuable
discussions.}
\begin{abstract}
  We approximate the homogenization of fully nonlinear, convex, uniformly elliptic Partial
  Differential Equations in the periodic setting, using a variational
  formula for the optimal invariant measure, which may be derived via
  Legendre-Fenchel duality. The variational formula expresses $\Hb$ as an
  average of the operator against the optimal invariant measure, generalizing
  the linear case. Several nontrivial analytic formulas for $\Hb$ are obtained.
  These formulas are compared to numerical simulations, using both PDE and
  variational methods. We also perform a numerical study of convergence rates
  for homogenization in the periodic and random setting and compare these to
  theoretical results.
\end{abstract}

\maketitle


\section{Introduction}
We consider homogenization of the periodic, convex, uniformly elliptic Hamilton-Jacobi-Bellman operator
\begin{equation}\label{eq:HJB}
  H(Q, y) = \sup_{\al \in \Ac}  L_\al(Q,y)   =\sup_{\al \in \Ac} \left\{
  -A(y,\al):Q + h(y,\al) \right\}.
\end{equation}
Note that $H(Q,y)$ is convex in $Q$. Let $\Ac$ be a convex compact control set, and let $A:
\Td \times \Ac \to \Sd$ and $h:\Td\times \Ac \to \R$ be continuous, where $\Sd$
is the space of $d\times d$ symmetric matrices, and $\Td$ is the $d$-dimensional
torus.  We assume that $A$ is uniformly elliptic, with $0 \ll  \lambda I  \ll A
\ll  \Lambda I$.

We will make use of the following result, which follows from
\Cref{th:main} below.  Consider an admissible control
$\alpha(y)$.  Let $L_\alpha$ be the corresponding linear operator, and let
$\overline {L_\alpha}(Q)$ the homogenized linear operator.  Then
\begin{equation}
  \label{LbarHbarEst}
  \overline {L_\alpha}(Q) \leq \Hb
\end{equation}
with equality when $\alpha(y)$ is optimal, or, equivalently, when the control  corresponds to the linearization about the corresponding solution of the cell problem, $Q + D^2u^Q(y)$.

In Section~\ref{sec:Analytics} we consider three example problems. One example is the maximum of two linear operators.  In this case, we obtain a formula for $\Hb$, which is new (as far as we know).  The second example is one dimensional, but with a quadratic nonlinearity.  In this case, by considering constant controls, we find a lower bound for $\Hb$ which is numerically verified to be sharp.

Finally, we consider a two dimensional Pucci operator on stripes.   In
\cite{ObermanFinlay_Homog}  we homogenized Pucci operators, mainly with
checkerboard coefficients.  There we did not require convexity of the operator.
We obtained accurate results for values of $Q$ away from the singularities of
the operators by simply linearizing the operator about $Q$. However, for stripes
coefficients, the linearization about $Q$ is not accurate.  Here, we linearize
about a control, and find the optimal constant control which corresponds to a
control direction which depends on both the eigenvectors of $Q$ and the
orientation of the stripes.  When compared numerically to  $\overline H(Q)$, this control gives very accurate results, away from the singularities.  Near the singularities, there is still a small nonzero error.    In \cite{ObermanFinlay_Homog} we also established upper bounds for the linear homogenization error.  These estimates included a term which decreased with the distance to the singular set of the operators.  Similar results apply here as well.

We compared our estimates for $\Hb$ to numerical results.
We computed $\Hb$ numerically using two methods: by solving the PDE for the
cell problem, and by using linear programming to solve for the invariant
measure.    We also considered the case of random coefficients, and we found that very similar
formulas for $\overline {L_\alpha}(Q)$ hold in the random setting.

We also computed rates of convergence for $\Hb$.  In the periodic case, we obtained second order convergence rates,
$\mathcal O(\e^2)$, in one dimension.  In the random setting we obtained a convergence rate of
$\mathcal O(\e^{1/2})$,  again in one dimension.  These are consistent with the theoretical results we mention below.

\subsection{Related work}
We know of few analytical solutions $\Hb$, other than the
formula for a separable Hamiltonian in one dimension which can be found in the
early paper~\cite{Lions1987}.  In \cite{ObermanFinlay_Homog} we obtained a formula for a separable second order operator in one dimension
\[
  \overline {a(y)F_0(Q)} = \HM(a)F_0(Q), \qquad y \in \R^1
\]
where $\HM(a)$ is the harmonic mean of the coefficient $a(y)$.  In general, the
harmonic mean formula on the right hand side is a lower bound.  However, in
\cite{ObermanFinlay_Homog} we found examples of Pucci type operators where the numerically compute value of $\Hb$ is very close to the right hand side, for values of $Q$ away from the corners of the operator.

For a general reference on theoretical and numerical homogenization in this
context, we refer to the review paper \cite{Engquist2008}.

A numerical method which uses the $\inf\sup$ formula for the first order case was developed in \cite{GO04}.
In~\cite{oberman2009homogenization} we
studied homogenization of convex (first order) Hamilton-Jacobi equations; some
exact formulas in the periodic setting can be found there.  Recently
\cite{Cacace2016} studied numerical homogenization of mainly first order
equations, along with one dimensional second order equations.  In
\cite{caffarelli2008numerical} the problem of homogenization of a Pucci type
equation with checkerboard coefficients was studied.  In that case, our results are close to,
but different from theirs, see \cite{ObermanFinlay_Homog}.

In the random setting, the first qualitative homogenization results  for fully
nonlinear uniformly elliptic operators were obtained
in~\cite{caffarelli2005homogenization}, followed by \cite{caffarelli2010rates}
which established a logarithmic estimate for convergence rates
in strongly mixing environments. Algebraic convergence estimates were
established in \cite{Armstrong2014}, where it was shown that in a uniformly
mixing environment,
\[
  \mathbb P \left[\norm{u^\varepsilon - \bar u }_\infty \geq C \varepsilon^\beta \right]
  \leq C \varepsilon^\beta,
\]
where $C$ and $\beta$ are constants that do not depend on $\varepsilon$.   In
the periodic case~\cite{camilli_rates_2009}, proved that the order of
convergence for the  cell problem is $\mathcal O(\varepsilon^2)$, when the HJB
operator does not depend on first order terms or the macroscopic scale.

\subsection{Background Theory}
In the periodic uniformly convex setting, the homogenized operator can be obtained by solving the cell problem.
\begin{definition}[Cell problem]
  Given $H$ as in  \eqref{eq:HJB}, for each $Q \in \Sd$, there is a unique value
  $\Hb$ and a periodic function $u^Q(y)$  which is a viscosity solution of the
  cell problem
  \begin{equation}\label{def:cell_problem}
    H(Q + D^2u^Q(y), y) = \Hb.
  \end{equation}
\end{definition}
Because the operator $H$ is uniformly elliptic, one can show that both the value
$\Hb$ and the solution $u^Q$ exist and are unique, and that
$u^Q \in C^2(\Td)$~\cite{Evans1989}.

In the linear case, we may use the Fredholm
Alternative to find the invariant measure, and the homogenized operator is then
obtained by averaging against the invariant measure, see for example
\cite{bensoussan2011asymptotic} and \cite{Froese2009}.  That is, under an integral compatibility condition, there is a unique
invariant probability measure, $\rho$, which solves $ D^2:(A(y)\rho) = 0$, and
the homogenized PDE operator is $\bar L(Q) = \bar A : Q + \bar h$   where
$\bar A = \int_\Td A \dd \rho$ and $\bar h = \int_\Td h \dd \rho$.

In the nonlinear case, the homogenized operator may still be found by averaging against the
optimal invariant measure~\cite{Gomes2005}, or~\cite{Ishii2016}.

\begin{definition}[Optimal invariant measure]
  Let $\PrAB$ denote the space of Borel  probability measures on $\Td \times \Ac$.
  For $\rho \in \PrAB$, we say $\rho$ is an \emph{invariant measure} if
  $L_0^*\rho(y,\al) = 0$ hold in the weak sense, by which we mean
  \begin{equation}\label{eq:dualpde}
    \int_{\Ac \times \Td} A(y,\alpha(y)):D^2\varphi(y)  \dd \rho(y, \al) = 0, \quad
    \forall \varphi \in C^2(\Td).
  \end{equation}
  Define
  \begin{equation}\label{eq:dual}
    \Bar H_{LP}(Q) =  \sup_{\rho \in \PrAB}
    \left\{ \int L_\al(Q,y) \dd  \rho(y,\al)  ~ \Big |~  L_0^* \rho = 0 \right\},
  \end{equation}
\end{definition}

The result follows from duality and convex analysis, in particular the theorem
of Fenchel-Rockafeller.  The fact that \eqref{eq:dual} and \eqref{def:cell_problem}
give the same value  is established in Theorem~\ref{th:maintheorem}.
\begin{theorem}\label{th:maintheorem}\label{th:main}
  Let $\Hb$ be defined by \eqref{def:cell_problem} and $ H_{LP}(Q)$ be defined by
  \eqref{eq:dual}.  Then
  \[
    \Hb = \Bar H_{LP}(Q).
  \]
\end{theorem}

\begin{remark}
  The formula above expresses an optimal invariant measure as a maximizer of the
  functional in \eqref{eq:dual}, and the homogenized operator as the average of
  $L_\alpha(Q,y)$ against an optimal invariant measure.

  Note that while the optimal invariant measure depends on $Q$, the set of
  invariant measures does not.  This allows us to sometimes find $\bar H(Q)$ for
  all $Q$ once the invariant measures are determined.
  While $\bar H(Q)$ does not depend on how we represent $H(Q,y)$ in
  \eqref{eq:HJB}, the set of invariant measure does. So a more concise
  representation of the operator can lead to a smaller set of invariant measures.
\end{remark}

\section{Estimates from Linearization}\label{sec:Analytics}
In this section we apply \eqref{LbarHbarEst}, by considering specific operators where we can obtain analytical values for the homogenized linear operator $\overline {L_\alpha}(Q)$.

\subsection{Pucci type operators on stripes}\label{sec:2d-P}
Consider the following Pucci type operator.
\begin{definition} \label{def:our_pucci}
  Let $y \in \T^2$ and $Q \in \mathcal{S}^2$.  Write  $\lp(Q),
  \lm(Q)$ for the maximum and minimum eigenvalues of $Q$, respectively.
  Given $A(y)\geq a(y) > 0$, define the convex Pucci type operator
  \begin{equation}\label{eq:OurPucci}
    F^{A,a} (Q,y) =
    a(y)\tr Q + b(y) \lp^+(Q)
  \end{equation}
  where $b(y) = A(y)-a(y)$ and $t^+ := \max \left\{t,0\right\}$
\end{definition}
\begin{remark}
The operator  \eqref{eq:OurPucci} can be rewritten
\begin{equation}\label{eq:Pucci2}
  F^{A,a} (Q,y) = a(y) \tr Q   +   b(y) \sup_{\Abs{\vec v}=1, \vec v = 0}
  \left\{ \vec v^\trp
  Q \vec v\right\}
\end{equation}
 When $Q$ is
negative definite, the operator is linear, and we obtain the harmonic mean.
  The level sets of this operator have corners on the negative axes and on the positive diagonal in the  $\lambda_1-\lambda_2$ plane. Elsewhere, the operator is linear in $\lambda_1$ and $\lambda_2$.
\end{remark}

For the rest of the discussion we restrict to $Q$ with at least one positive
eigenvalue.

\begin{formula}[Pucci with stripes]\label{formula:stripes}
  Consider $F^{A,a} (Q,y)$ given by \eqref{eq:OurPucci} in dimension $d = 2$.
  Consider piecewise constant stripes, with
  \begin{equation}\label{stripes2}
    a(y_1) = 1,
    \qquad
    b(y_1) = \begin{cases}
      0,\quad 0 \leq y_1 \leq \frac{1}{2} \\
      b_0,\quad \frac{1}{2} \leq y_1 \leq 1
    \end{cases}
  \end{equation}
For $Q \not\preceq 0$ and for $t \in [0,1]$ define
\[
\overline {L_t}(Q) = \tr Q +   \frac{b_0}{2 + b_0 t}
    \left (
    q_{22} +t(q_{11} - q_{22}) + 2 q_{12} \sqrt{t - t^2}
    \right ),
\]
Then
  \begin{equation}\label{Hbarstripes2}
    \overline {F^{A,a}} (Q)  \geq  \sup_{t \in [0,1]} \overline {L_t}(Q),
        \quad Q \not\preceq 0
  \end{equation}
and
\[
   \HM (a) \tr Q,        \quad Q \preceq 0.
\]
\end{formula}

\begin{remark}
This formula is obtained by considering constant controls in the direction of a given unit vector.  We homogenize the linearization about this control, and then optimize over the choice of unit vector.  In two dimensions, these unit vectors are determined by a single parameter.

  The term in \eqref{Hbarstripes2} can be simplified further analytically, but
  the formula becomes complicated.  It is more convenient to solve it
  numerically using one variable equation solvers.
\end{remark}

\begin{proof}
We need only consider the case $Q \not \preceq 0$ since the operator is linear otherwise.

1.  Since we are on stripes, we restrict to invariant measures which depend only on $y_1$.    Given any choice of control
  $\alpha(y_1)$, the equation for the invariant measure is
  \[
    \partial_{11} \left( c(y_1) p_{\alpha(y_1)}(y_1) \right ) = 0,
    \quad  c(y_1) = a(y_1) + b(y_1) \alpha(y_1).
  \]
  and the corresponding invariant measure is
  $
  p_{\alpha(\cdot)}(y_1) = \frac{\HM(c(y_1))}{c(y_1)}.
  $

2. Next with coefficients given by \eqref{stripes2}, restrict to $\alpha(y_1)$
  constant on $b(y_1) = b_0$.

3.   In two dimensions,  parameterize the unit vectors $(\cos \theta, \sin \theta)$
  by $ \vec v_\alpha = (\sqrt \alpha , \sqrt{1-\alpha})$, for $\alpha = \cos^2
  \theta$, and write
  \begin{equation}\label{Balpha}
    B_\alpha =  \vec v_\alpha \vec v_\alpha^\trp =
    \left [
      \begin{matrix}
        \alpha & \sqrt{\alpha(1-\alpha)} \\
        \sqrt{\alpha(1-\alpha)} & 1-\alpha
      \end{matrix}
    \right ]
  \end{equation}
  so that \eqref{eq:Pucci2} becomes
  \[
    F^{A,a} (Q,y) = \tr Q  + b(y) \sup_{ \alpha \in [0,1] } \tr B_\alpha  Q
  \]
  For these measures, the homogenized linear operator becomes
  \begin{equation}\label{barF}
    \HM(a) \tr Q +  \int_{\T^1} b(y_1) \tr \left ( B_{\alpha}  Q \right) p_{\alpha}(y)	dy
  \end{equation}
  after integrating out the invariant measure from the Laplacian term.

  Use the representation $B_\alpha$ given by~\eqref{Balpha} to simplify
  \eqref{barF} to obtain the expression which is maximized in
  \eqref{Hbarstripes2}.
\end{proof}

\subsection{Maximum of two linear operators}

\begin{definition}\label{defn:HmaxLin}
  Given  a (constant) symmetric positive definite matrix, $A$, positive
  functions $a_0(y), a_1(y) > 0$, and the constant $h$. Define
  \begin{equation}
    H(Q,y) = \max \{ a_0(y), a_0(y) + a_1(y) \} A:Q + h
    \label{eq:max2lin}
  \end{equation}
\end{definition}

\begin{formula}[Maximum of two linear operators]\label{formula:3}
  Let $H(Q,y)$ be given by \eqref{eq:max2lin}.  Then
  \begin{equation}\label{max2}
    \bar H(Q) = \max \{  -\HM(a_0) A:Q, -\HM(a_0 + a_1) A:Q\} + h
  \end{equation}
\end{formula}

\begin{proof}
  Write
  \[
    H(Q,y) = \max_{\alpha \in [0,1]}
 L_\alpha(Q,y),
    \quad
  L_\alpha(Q,y)\equiv    -(a_0(y) + \alpha a_1(y)) A:Q + h
    \]
1. For any  choice $\alpha(y)$, the corresponding invariant measure defined by
  $L_{\alpha(y)}(Q,y)$, is given by
  \[
    p_\alpha(y) = \frac{\HM(b(y, \alpha(y))}{b(y, \alpha(y))},
    \quad b(y, \alpha(y)) = a_0(y) + \alpha(y)a_1(y).
  \]
and
  \[
  \overline {L_{\alpha(y)}}(Q) =
   \int_\Td L_{\alpha(y)}(Q,y)\dd p_\alpha(y)
    =\HM(b(y, \alpha(y))) A:Q + h
  \]
So from \eqref{LbarHbarEst}, we have
\[
\overline H (Q) = \sup_{\alpha(\cdot)}   \overline {L_{\alpha(y)}}(Q)
\]
Notice that  $b = (a_0(y) + \alpha a_1(y))$ is increasing in $\al$ for each $y$.
Moreover, it is easy to verify that the Harmonic Mean is an increasing function of $b$.  Thus, depending on the sign of $A:Q$, the optimal value  is achieved by either $\al(y) \equiv 0$ or $\al(y) \equiv 1$, accordingly.   This gives~\eqref{max2}.
\end{proof}

\subsection{A one dimensional quadratic operator}\label{sec:quad}

\begin{definition}\label{Quadratic operator}
  Consider for constants, $c$, $a > 0$ and for the function $b(y)\geq 0$,
  \[
    H(Q,y) = aQ + b(y)(Q^+)^2 - c.
  \]
\end{definition}
It easy to verify that
  \[
    H(Q,y) = \max_{\alpha\in [0,Q]} L_\alpha(Q,y), \quad     L_\alpha(Q,y)  =  \big(a + 2b(y)\alpha\big) Q - \big(b(y) \alpha^2
    +c\big).
  \]

\begin{formula}\label{formula:quadratic}
  Let $H(Q,y)$ be given as in Definition~\ref{Quadratic operator}. Consider the one dimensional case, $d=1$, and suppose $b$ is piecewise constant,
  \begin{equation}
    b(y) =
    \begin{cases}
      0,\quad 0 \leq y \leq \frac{1}{2} \\
      b_0,\quad \frac{1}{2} \leq y \leq 1
    \end{cases}
    \label{eq:1d-Q-b}
  \end{equation}
  Then
  \begin{equation}\label{eq:1d-quad-barH}
    \Hb \geq  a(Q+Q^+) - c + \frac{a^2}{b_0} -
    \frac{1}{b_0}\sqrt{a^3(a+2b_0Q^+)}.
  \end{equation}
\end{formula}

\begin{proof}
Consider constant controls $\alpha(y) \equiv \alpha$. In this case the invariant measure
  $p_{\alpha}(y)$ is given by
  \begin{equation}
    p_{\alpha}(y) = \begin{dcases}
      \frac{a+2 b_0 \alpha}{a+b_0 \alpha},\quad 0 \leq y \leq \frac{1}{2} \\
      \frac{a}{a+b_0 \alpha},\quad \frac{1}{2} \leq y \leq 1.
    \end{dcases}
    \label{eq:1d-quad-measure}
  \end{equation}
and
  \begin{align}
    \bar L_\alpha(Q) &=  \langle L_\alpha, p_\alpha \rangle\\
    &=  \frac{a \left( a+2b_0\alpha \right)}{a+b_0\alpha}
    \left[ Q - \frac{1}{2}\left( \frac{c}{a} +
    \frac{b_0 \alpha^2 + c}{a + 2 b_0 \alpha} \right) \right].
    \label{eq:1d-Q-ugly}
  \end{align}

By \eqref{LbarHbarEst},
\[
\Hb \geq \max_{\alpha\in [0,Q]} \overline {L_\alpha}(Q)
\]
Next, maximize over $\alpha$.  This is accomplished by solving for the roots of the derivative of this expression with respect to $\alpha$.  We obtain
  \begin{equation}
    \alpha^*(Q) = \frac{1}{b_0}\left(-a + \sqrt{a(a+2b_0 Q^+)}\right).
    \label{eq:1d-quad-alpha}
  \end{equation}

Thus the estimate is givenby \eqref{eq:1d-Q-ugly} with $\alpha$ given by \eqref{eq:1d-quad-alpha}.  Upon simplification, we obtain~\eqref{eq:1d-quad-barH}.
\end{proof}

\section{Numerical results}
Here we compare the results of \Cref{sec:Analytics} with the numerical
homogenization of the operators.

\begin{remark}[Numerical methods]
  $\Hb$ was computed with two methods. In the first, the equation \eqref{eq:HJB}
  was discretized with finite differences. A steady state solution was
  computed iteratively by Euler step to the parabolic equation $u_t
  +H(Q+D^2u,y)$. We used a filtered scheme \cite{froese2013convergent} to choose
  between a monotone finite difference scheme and standard accurate finite differences.
  However the standard finite difference scheme was always chosen by the
  filtered scheme, likely because solutions are $C^2$ and periodic.

  We also computed $\bar H_{LP}(Q)$ by discretizing the control space $\Ac$ and
  formulating the problem \eqref{eq:dual}  as a discrete linear programming
  problem.
  Derivatives were discretized via standard second order finite differences. We
  then solved this LP using the package CVX \cite{cvx,gb08} with the SeDuMi
  solver \cite{sturm1999}.

  Throughout we set $h(y,\alpha) = c = 1$, and subtracted this constant from
  $\Hb$, so as to avoid trivial solutions.
\end{remark}

\subsection{Pucci type operator on stripes}

\begin{figure}[t]
  \centering
  \begin{subfigure}[b]{0.48\textwidth}
    \includegraphics[width=\textwidth]{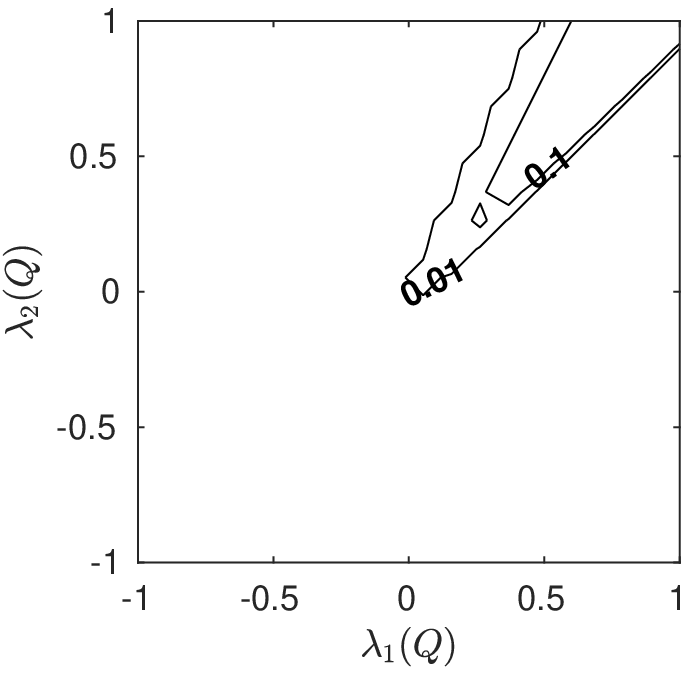}
    \caption{Optimal Linearization error}
    \label{fig:stripes-lp}
  \end{subfigure}
  \hfill
  \begin{subfigure}[b]{0.48\textwidth}
    \includegraphics[width=\textwidth]{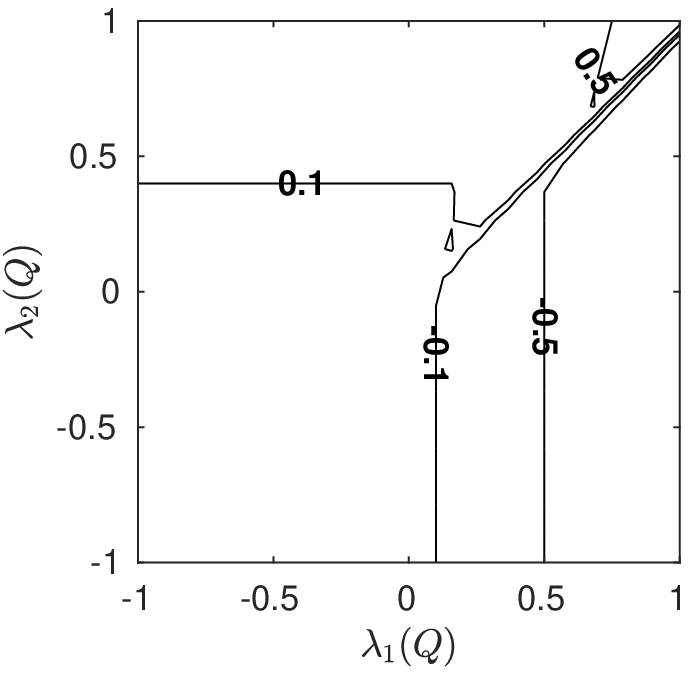}
    \caption{Linearization error}
    \label{fig:stripes-lin}
  \end{subfigure}
  \label{fig:stripes}
  \caption{Comparison between homogenization error using an invariant measure,
    and by homogenizing the linearized operator. In this example the operator is given by
    \eqref{stripes2}, with $b_0=2$. Here $Q = \text{diag}
    (\lambda_1,\lambda_2)$.  In the third quadrant, the operator is linear, and the error was zero up to machine precision.
    Figure \ref{fig:stripes-lp}: error of Formula \ref{formula:stripes}, the
    error is \num{1e-8} is most of the domain, with the \num{1e-2} level set shown.
        Figure \ref{fig:stripes-lin}: error of homogenizing the linearized
    operator.  There the error is order one, outside the third quadrant.
    }
    \end{figure}
We compared the analytical formula, Homogenization Formula \ref{formula:stripes},
with numerically homogenized values. This required solving a one-variable
optimization problem. We used piecewise constant
coefficients, where the operator was  either $\tr Q$, or $F^{3,1}(Q) = 1 \tr Q +
2 \lp^+(Q)$. The error profile of the analytic lower bound against
the numerical homogenization is plotted in \Cref{fig:stripes-lp}, for a set of
diagonal $Q$. In the vicinity
of the line $\lambda_1(Q) = \lambda_2(Q)$ the error is on the order of \num{1e-1};
elsewhere the error is less than \num{1e-2}.

We contrast this homogenization approach with the method of homogenizing the
linearized operator \cite{ObermanFinlay_Homog}. The homogenizing error
by first linearizing the operator is much greater than the error given by
Formula \ref{formula:stripes}, as can be seen by comparing
\Cref{fig:stripes-lp,fig:stripes-lin}.

There is a symmetry in $\bar H(Q)$.  We represent
\[
  Q = R_\phi^T  \text{diag} (\lambda_1,\lambda_2) R_\phi
\]
where $R_\phi$ is a rotation matrix.  When $\phi = \pi/4$,  the orientation of
the stripes is at an equal angle to the eigenvectors of $Q$, then $\bar H(Q)$ is
symmetric about $\lambda_1 = \lambda_2$.   More generally $\bar H(Q)$ is
symmetric under reflections in the angle about the same line of symmetry:
\[
  \bar H( Q |_{\phi = \pi/4 - \gamma} )  = \bar H( Q|_{\pi/4 + \gamma} ),
  \qquad
  \text{for $\abs{\gamma} \leq \pi/4$.}
\]

\subsection{Maximum of two linear operators, in one and two dimensions}\label{sec:1d-P}
\begin{figure}[t]
  \centering
  \includegraphics[width=.5\textwidth]{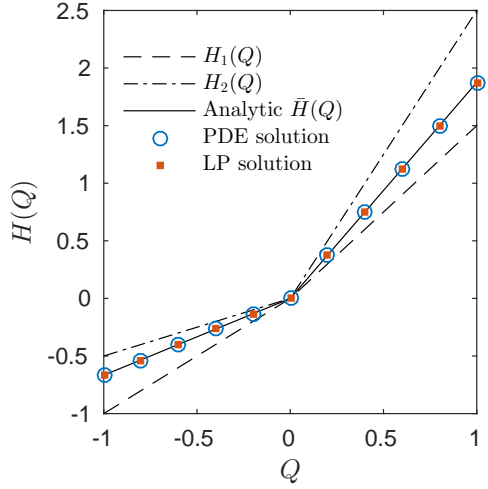}
  \caption{Validation of Formula \ref{formula:3}, homogenization of the maximum of two
  linear operators, \eqref{eq:max2lin}.
  Lines represent $\bar H(Q)$ and each of the operators $H_i(Q)$.
}
\label{fig:p-1d}
\end{figure}

\begin{figure}[t]
  \centering
    \includegraphics[width=0.5\textwidth]{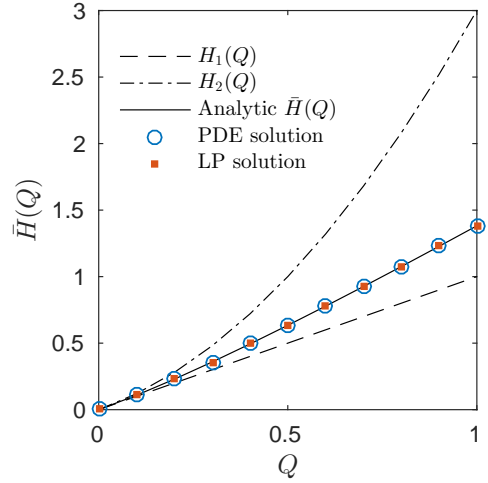}
  \caption{Validation of Formula~\ref{formula:quadratic}.
 Value of $H_1(Q), H_2(Q)$, the
  numerically homogenized operators, and the analytic homogenized operator.
  }
  \label{fig:1d-quad-ops}
\end{figure}

We numerically validated Homogenization Formula~\ref{formula:3}, for the maximum
of two linear operators.  We considered the case when the dimension is one, and
we took $A = 1$.  The interval $[0,1]$ was discretized into $20$ equal sized
pieces. The
coefficients $a_0(y)$ and $a_1(y)$ were piecewise constant on these equal-sized
pieces.  We took $h(y)$ to be constant.  In Figure \ref{fig:p-1d} we let $a_0$
alternate between $1$ and $\frac{1}{2}$, and let $a_1$ alternate between
$\frac{3}{2}$ and  $\frac{5}{2}$.   The values of the analytically homogenized
operator are indistinguishable from the numerically computed values, for
discrete values of $Q$, using both the direct method and the dual method.  Even
at the discontinuity $Q=0$, the formula agrees with the numerical homogenization
up to machine precision.  For reference, we also plotted $H_1(Q) = \min_y
H(Q,y)$ and $H_2(Q) = \max_y H(Q,y)$.  We computed many different examples and
obtained similar results.  (Note in this example, the invariant measure is
piecewise quadratic, so the numerical method is very accurate.)  We also
visualized the numerical invariant measure, and found that it agreed with our
formula.

We also numerically validated Formula \ref{formula:3} in two dimensions, and
obtained similar results: in this case the analytic formula and the numerical simulations
agree up to \num{1e-12}.

\subsection{The quadratic operator}
Next we considered the example from \S\ref{sec:quad},  Homogenization
Formula~\ref{formula:quadratic}, for the operator
\[
  H(Q,y) = aQ + b(y)(Q^+)^2 - c.
\]
Here $c$ is a constant. We numerically homogenized this operator on the periodic
domain $[0,1]$, divided into 20 pieces.  The coefficients are piecewise constant
on equal intervals.  As illustrated in Figure \ref{fig:1d-quad-ops} the analytic
homogenization and the numerically homogenized operator are indistinguishable.
As in the previous operator (maximum of two linear operators), even at $Q=0$ the formula agrees with the numerical
homogenization up to machine precision.
Again, we also plot $H_1(Q) = \min_y H(Q,y) = a Q$ and $H_2(Q) = \max_y
H(Q,y) = aQ + b_0(Q^+)^2$.

\section{Numerical rates of convergence in the periodic and random case}
Using the exact analytical formulas of Section
\ref{sec:Analytics}~(Formula \ref{formula:3}  and Formula
\ref{formula:quadratic}), we investigate empirical rates of convergence of the
small-scale solutions $u^\e$ to the solution $\bar u$ of the homogenized
operator.   Although our theoretical results were for the periodic case, we
found that the same formulas applied in the random case.  This allows us to
study empirical convergence rates in the random case as well.

We solved the Dirichlet problem with zero boundary conditions on the interval
$[0,1]$
for the two different operators, in both the random and periodic case.
The operators were the maximum of two linear operators,
Formula~\ref{formula:3}; and the quadratic operator,
Formula~\ref{formula:quadratic}.  These are both one dimensional examples.

We used a sequence of decreasing cell sizes, of width $\e$.   We used 1 grid
point per cell.

We also solved the same problem with the homogenized
operator.  Numerically we obtained two solutions, $u^\e$,  and
$\bar u$ corresponding to
\begin{equation}
\begin{dcases}
H^\e(D^2 u^\e(x),x) = 1 \\
u^\e(x) = 0, \quad x \in \partial[0,1]
\end{dcases}
\qquad\text{ and }
\quad\begin{dcases}
\bar H(D^2\bar u(x),x) = 1 \\
\bar u(x) = 0, \quad x \in \partial[0,1].
\end{dcases}
\end{equation}

We chose coefficients which were piecewise constant.  Let $H^\e(Q,x)$ be the
operator parameterized by checkerboard square width $\e$. We checked convergence
for both the periodic case, and the random case. In the periodic case, the
unhomogenized operator alternates between two constituent operators $H_1$ and
$H_2$ between the checkerboard cells. In the random case, in each checkerboard
square we randomly sample from the two constituent operators with probability
$\frac{1}{2}$.

In the random case, we observed convergence rates consistent with $\mathcal O(\e
^\frac{1}{2})$ in the sup-norm for both operators.  In the periodic case, we
observed convergence rates consistent with $\mathcal O(\e^2)$.

\begin{figure}[t]
  \begin{subfigure}[b]{0.36\textwidth}
    \includegraphics[width=\textwidth]{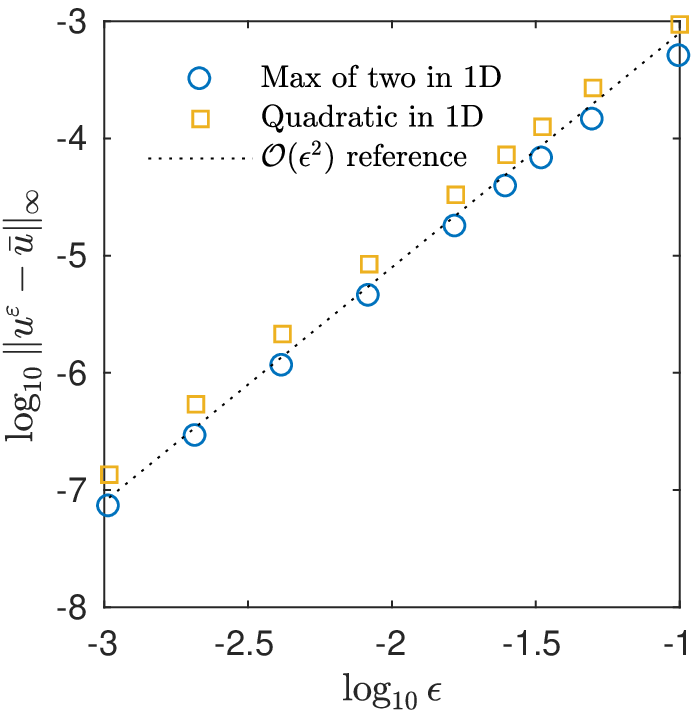}
    \caption{Periodic convergence}
    \label{fig:p-convg}
  \end{subfigure}
  \begin{subfigure}[b]{0.61\textwidth}
    \includegraphics[width=\textwidth]{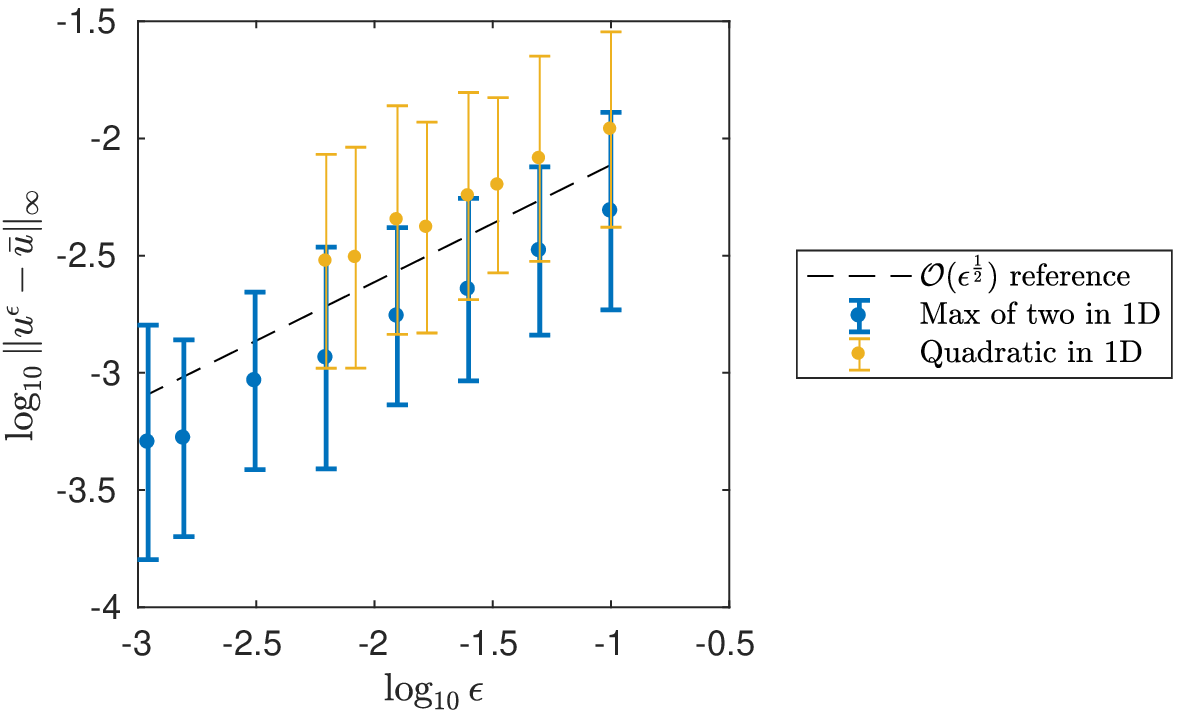}
    \caption{Random convergence}
    \label{fig:rand-convg}
  \end{subfigure}
  \caption{Figure \ref{fig:convg}: Periodic coefficients: rate of convergence
  $u^\e \to \bar u$.  Figure \ref{fig:rand-convg}: Random coefficients:  rate
  of convergence.  We plot 90\% confidence intervals for a normal distribution.}
\label{fig:convg}
\end{figure}

Figure \ref{fig:convg} presents the observed rates of convergence as $u^\e \to
\bar u$ in the sup-norm.
In the periodic setting, the order is nearly $\mathcal O(\e^2)$:
we estimate that the order of convergence is $\mathcal
O(\e^{1.95})$. In the random setting, we solved
each problem 20 times, drawing the random checkerboard anew at each iteration.
We then used least squares to estimate the order of convergence. We
summarize these convergence estimates in Table \ref{tab:poly-ord}. It appears
that convergence in the sup-norm is roughly $\mathcal
O(\e^\frac{1}{2})$ on the random checkerboard. We also measured the errors in the $\ell^2$ and
$\ell^1$ norms.
\begin{table}
\centering
\begin{tabular}{l | p{1.5cm} | p{1.5cm} |  p{1.5cm} | p{1.5cm} }
\hline
Operator & Periodic,\newline $\norm{\cdot}_\infty$ & Random, \newline$\norm{\cdot}_\infty$ &
Random, \newline$\norm{\cdot}_2$ & Random, \newline$\norm{\cdot}_1$\\
\hline
Max of two in 1D & 1.95 & 0.51 & 0.49 & 0.50 \\
Quadratic in 1D & 1.95 & 0.48 & 0.42 & 0.41 \\
\hline
\end{tabular}
\caption{Empirical rates of convergence for the two operators.}
\label{tab:poly-ord}
\end{table}

\section{Conclusions}
In this article we investigated the accuracy of approximating nonlinear homogenization by the homogenization of a linearization of the operator.  In previous work \cite{ObermanFinlay_Homog}, we simply linearized about a constant.  There, we obtained very accurate for checkerboard type coefficients, but significant errors in the case of stripes.   In this article, we restricted to convex operators.  This allowed us to write operators as the supremum of linear operators.  For any linearization over a choice of control $\alpha(y)$
\[
  \overline {L_\alpha}(Q) \leq \Hb
\]
with equality when $\alpha(y)$ is optimal.

We applied this formula to three examples.  For the example of a maximum of two
linear operators, we obtained an exact result, given by the maximum of two
harmonic means (see~\eqref{max2}).  For a quadratically nonlinear one
dimensional operator, we restricted to piecewise constant controls and optimized
over the value of the control.  This results in a lower bound which was verified
by numerical simulations to be exact.

Finally, we consider the Pucci-type operator with stripe coefficients.  In this case, the controls depended on a choice of direction vector, which in two dimensions resulted in a one parameter optimization problem for $\overline {L_\alpha}(Q)$.   The solution of this problem
was verified by numerical simulations to be nearly exact over parameter values away from the singularities of the operator.  For other values of $Q$ it achieved a small (a few percentages) relative error.

We also consider the numerical convergence rates of the homogenization problem in the scale parameter, obtaining results consistent with recent analytical results, in both the periodic and random case.

\bibliographystyle{alpha}
\bibliography{HomogLP}

\begin{thebibliography}{CSW05}

\bibitem[AS14]{Armstrong2014}
Scott~N. Armstrong and Charles~K. Smart.
\newblock Quantitative {Stochastic} {Homogenization} of {Elliptic} {Equations}
  in {Nondivergence} {Form}.
\newblock {\em Archive for Rational Mechanics and Analysis}, 214(3):867--911,
  December 2014.

\bibitem[BLP11]{bensoussan2011asymptotic}
Alain Bensoussan, Jacques-Louis Lions, and George Papanicolaou.
\newblock {\em Asymptotic analysis for periodic structures}, volume 374.
\newblock American Mathematical Soc., 2011.

\bibitem[CC16]{Cacace2016}
Simone Cacace and Fabio Camilli.
\newblock Ergodic problems for {Hamilton}-{Jacobi} equations: yet another but
  efficient numerical method.
\newblock {\em arXiv preprint arXiv:1601.07107}, 2016.

\bibitem[CG08]{caffarelli2008numerical}
LA~Caffarelli and Roland Glowinski.
\newblock Numerical solution of the dirichlet problem for a pucci equation in
  dimension two. application to homogenization.
\newblock {\em Journal of Numerical Mathematics}, 16(3):185--216, 2008.

\bibitem[CM09]{camilli_rates_2009}
Fabio Camilli and Claudio Marchi.
\newblock Rates of convergence in periodic homogenization of fully nonlinear
  uniformly elliptic {PDEs}.
\newblock {\em Nonlinearity}, 22(6):1481, 2009.

\bibitem[CS10]{caffarelli2010rates}
Luis~A Caffarelli and Panagiotis~E Souganidis.
\newblock Rates of convergence for the homogenization of fully nonlinear
  uniformly elliptic pde in random media.
\newblock {\em Inventiones mathematicae}, 180(2):301--360, 2010.

\bibitem[CSW05]{caffarelli2005homogenization}
Luis~A Caffarelli, Panagiotis~E Souganidis, and Lihe Wang.
\newblock Homogenization of fully nonlinear, uniformly elliptic and parabolic
  partial differential equations in stationary ergodic media.
\newblock {\em Communications on pure and applied mathematics}, 58(3):319--361,
  2005.

\bibitem[ES08]{Engquist2008}
Bj{\"o}rn Engquist and Panagiotis~E Souganidis.
\newblock Asymptotic and numerical homogenization.
\newblock {\em Acta Numerica}, 17:147--190, 2008.

\bibitem[Eva89]{Evans1989}
Lawrence~C Evans.
\newblock The perturbed test function method for viscosity solutions of
  nonlinear {PDE}.
\newblock {\em Proceedings of the Royal Society of Edinburgh: Section A
  Mathematics}, 111(3-4):359--375, 1989.

\bibitem[FO09]{Froese2009}
Brittany~D. Froese and Adam~M. Oberman.
\newblock Numerical averaging of non-divergence structure elliptic operators.
\newblock {\em Communications in Mathematical Sciences}, 7(4):785--804, 2009.

\bibitem[FO13]{froese2013convergent}
Brittany~D. Froese and Adam~M. Oberman.
\newblock Convergent filtered schemes for the {M}onge-{A}mp\`ere partial
  differential equation.
\newblock {\em SIAM J. Numer. Anal.}, 51(1):423--444, 2013.

\bibitem[FO]{ObermanFinlay_Homog}
Chris Finlay and Adam~M Oberman.
\newblock Approximate homogenization of fully nonlinear elliptic {PDE}s:
estimates and numerical results for Pucci type equations.

\bibitem[GB08]{gb08}
Michael Grant and Stephen Boyd.
\newblock Graph implementations for nonsmooth convex programs.
\newblock In V.~Blondel, S.~Boyd, and H.~Kimura, editors, {\em Recent Advances
  in Learning and Control}, Lecture Notes in Control and Information Sciences,
  pages 95--110. Springer-Verlag Limited, 2008.
\newblock \url{http://stanford.edu/~boyd/graph_dcp.html}.

\bibitem[GB14]{cvx}
Michael Grant and Stephen Boyd.
\newblock {CVX}: Matlab software for disciplined convex programming, version
  2.1.
\newblock \url{http://cvxr.com/cvx}, March 2014.

\bibitem[GO04]{GO04}
Diogo~A Gomes and Adam~M Oberman.
\newblock Computing the effective {Hamiltonian} using a variational approach.
\newblock {\em SIAM journal on control and optimization}, 43(3):792--812, 2004.

\bibitem[Gom05]{Gomes2005}
Diogo~Aguiar Gomes.
\newblock {\em Trends in Partial Differential Equations of Mathematical
  Physics}, chapter Duality Principles for Fully Nonlinear Elliptic Equations,
  pages 125--136.
\newblock Birkh{\"a}user Basel, Basel, 2005.

\bibitem[IMT16]{Ishii2016}
Hitoshi Ishii, Hiroyoshi Mitake, and Hung~V. Tran.
\newblock The vanishing discount problem and viscosity mather measures. part 1:
  The problem on a torus.
\newblock {\em Journal de Math{\'e}matiques Pures et Appliqu{\'e}es}, 2016.

\bibitem[LPV87]{Lions1987}
Pierre-Louis Lions, George Papanicolaou, and Srinivasa~RS Varadhan.
\newblock Homogenization of {H}amilton-{J}acobi equations.
\newblock 1987.

\bibitem[OTV09]{oberman2009homogenization}
Adam~M Oberman, Ryo Takei, and Alexander Vladimirsky.
\newblock Homogenization of metric {H}amilton-{J}acobi equations.
\newblock {\em Multiscale Modeling \& Simulation}, 8(1):269--295, 2009.

\bibitem[Stu99]{sturm1999}
Jos~F Sturm.
\newblock Using {SeDuMi} 1.02, a {MATLAB} toolbox for optimization over
  symmetric cones.
\newblock {\em Optimization methods and software}, 11(1-4):625--653, 1999.

\end{thebibliography}
\end{document}